\documentclass{article}

\usepackage{amsmath}
\usepackage{amssymb}

\newtheorem{thm}{Theorem}

\title{Shortest Distance in Modular Hyperbola and Least Quadratic Nonresidue}
\author{Tsz Ho Chan}
\date{}

\begin{document}
\maketitle

\begin{abstract}
In this paper, we study how small a box contains at least two points from a modular hyperbola $x y \equiv c \pmod p$. There are two such points in a square of side length $p^{1/4 + \epsilon}$. Furthermore, it turns out that either there are two such points in a square of side length $p^{1/6 + \epsilon}$ or the least quadratic nonresidue is less than $p^{1/(6 \sqrt{e}) + \epsilon}$.
\end{abstract}

\section{Introduction and Main results}
Let $p > 2$ be a prime and $(c, p) = 1$. We consider the modular hyperbola
\[
H_c := \{ (x, y): x y \equiv c \pmod p \}.
\]
We are interested in the shortest distance between two points in $H_c$. Rather than distances, we consider how small a box
\[
B(X, Y; H) := \{ (x, y): X + 1 \le x \le X + H \pmod p, Y + 1 \le y \le Y + H \pmod p \}
\]
contains two points in $H_c$ where $X$ and $Y$ run over $0$, $1$, ..., $p-1$.

By H\"{o}lder's inequality and Weil's bound on character sum, we have
\begin{thm} \label{thm1}
For all $(c, p) = 1$, there exist some $0 \le X, Y \le p-1$ so that
\[
| H_c \cap B(X, Y; H) | \ge 2
\]
if $H \gg_\epsilon p^{1/4 + \epsilon}$.
\end{thm}
Now let us switch to the subject of the least quadratic nonresidue modulo $p$. Many people have been interested in the upper bound for $n_p$. By Polya-Vinogradov bound on character sums, we have
\[
n_p \ll p^{1/2} \log p.
\]
Vinogradov \cite{V} applied a trick and got
\[
n_p \ll_\epsilon p^{1/(2 \sqrt{e}) + \epsilon}.
\]
Burgess \cite{B} proved a new bound on short character sum which together with Vinogradov's trick yielded
\[
n_p \ll_\epsilon p^{1/(4 \sqrt{e}) + \epsilon}.
\]

Recall a recent result of Heath-Brown \cite{H} and Shao \cite{S} on mean-value estimates of character sums:
\begin{thm} \label{Shao}
Given $H \le p$, a positive integer and any $\epsilon > 0$. Suppose that $0 \le N_1 < N_2 < ... < N_J < p$ are integers satisfying $N_{j+1} - N_j \ge H$ for $1 \le j < J$. Then
\[
\sum_{j = 1}^{J} \max_{h \le H} |S(N_j; h)|^{2r} \ll_{\epsilon, r} H^{2r-2} p^{1/2 + 1/(2r) + \epsilon}
\]
where
\[
S(N; H) := \sum_{N < n \le N + H} \chi(n).
\]
and $\chi$ is any non-principal character modulo $p$.
\end{thm}

Applying the above theorem, we can show that
\begin{thm} \label{thm2}
For any $\epsilon  > 0$, we have either
\[
n_p \ll_\epsilon p^{1 / (6 \sqrt{e}) + \epsilon}
\]
or, for any $(c, p) = 1$ and $H \gg_\epsilon p^{1/6 + \epsilon}$,
\[
| H_c \cap B(X, Y; H) | \ge 2
\]
for some $0 \le X, Y \le p-1$.
\end{thm}

It is probably the case that the above two statements are true simultaneously. The paper is organized as follows. In section 2, we give a basic argument transforming the existence of two close points in the modular hyperbola to a certain equality in Legendre symbol. Then we prove Theorem \ref{thm1} in section 3 and Theorem \ref{thm2} in section 4.

\bigskip

{\bf Some Notations} Throughout the paper, $p$ stands for a prime. The symbol $|S|$ denotes the number of elements in the set $S$. We also use the Legendre symbol $(\frac{\cdot}{p})$. The notation $f(x) = o(g(x))$ means that the ratio $f(x) / g(x)$ is going to zero as $x, p \rightarrow \infty$. The notations $f(x) \ll g(x)$, $g(x) \gg f(x)$ and $f(x) = O(g(x))$ are equivalent to $|f(x)| \leq C g(x)$ for some constant $C > 0$. Finally, $f(x) \ll_{\lambda_1, ..., \lambda_k} g(x)$, $g(x) \gg_{\lambda_1, ..., \lambda_k} f(x)$ and $f(x) = O_{\lambda_1, ..., \lambda_k} (g(x))$ mean that the implicit constant $C$ may depend on $\lambda_1$, ..., $\lambda_k$.

\section{The Basic Argument} \label{basic}
For $(c, p) = 1$, suppose $| H_c \cap B(X, Y; H) | \ge 2$ for some $0 \le X, Y \le p-1$. This means that
\begin{equation} \label{start}
x y \equiv c \pmod p, \text{ and } (x + a) (y + b) \equiv c \pmod p
\end{equation}
for some $1 \le x, y \le p-1$ and $1 \le a, b \le H$. After some algebra, one can show that (\ref{start}) is equivalent to
\[
b x + a c \overline{x} + a b \equiv 0 \pmod p
\]
where $\overline{x}$ stands for the multiplicative inverse of $x$ modulo $p$ (i.e. $x \overline{x} \equiv 1 \pmod p$. This, in turn, is equivalent to
\[
(2 b x + a b)^2 \equiv (a b)^2 - 4 a b c \pmod p.
\]
Therefore $| H_c \cap B(X, Y; H) | \ge 2$ if and only if
\[
\Bigl( \frac{a b}{p} \Bigr) \Bigl( \frac{a b - 4 c}{p} \Bigr) = 1
\]
for some $1 \le a, b \le H$. We are going to restrict our attention to even $a = 2a'$'s and $b = 2b'$'s. So we want
\begin{equation} \label{conclude}
\Bigl( \frac{a' b'}{p} \Bigr) \Bigl( \frac{a' b' - c}{p} \Bigr) = 1 \text{ for some } 1 \le a', b' \le H/2.
\end{equation}

\section{Proof of Theorem \ref{thm1}}
Throughout this section, we assume that $H \gg_\epsilon p^{1/4 + \epsilon}$. We want to show that
\[
\sum_{a' \le H/2} \sum_{b' \le H/2} \Bigl( \frac{a' b'}{p} \Bigr) \Bigl( \frac{a' b' - c}{p} \Bigr) = o(H^2).
\]
Then either we have two pairs with $(\frac{a_1' b_1' - c}{p}) = 0 = (\frac{a_2' b_2' - c}{p})$ which gives Theorem \ref{thm1} automatically; or at most one such pair equal to $0$ which would imply that $( \frac{a_1' b_1'}{p} ) ( \frac{a_1' b_1' - c}{p} ) = 1$ and $( \frac{a_2' b_2'}{p} ) ( \frac{a_2' b_2' - c}{p} ) = -1$ for some $1 \le a_1', b_1', a_2', b_2' \le H/2$. However, we are going to restrict $a'$ to a special form, namely $a' = u v$ with $1 \le u \le U$, $1 \le v \le V$ and $U V = H/2$. So let us consider
\[
S := \sum_{u \le U} \sum_{v \le V} \sum_{b' \le H/2} \Bigl( \frac{u v b'}{p} \Bigr) \Bigl( \frac{u v b' - c}{p} \Bigr).
\]
Then
\[
S \le \sum_{v \le V} \sum_{b' \le H} \Big| \sum_{u \le U} \Bigl( \frac{u}{p} \Bigr) \Bigl( \frac{u v b' - c}{p} \Bigr) \Big| \ll_\epsilon p^\epsilon \sum_{n \le V H} \Big| \sum_{u \le U} \Bigl( \frac{u}{p} \Bigr) \Bigl( \frac{u n - c}{p} \Bigr) \Big|
\]
as the number of divisors of $n$ is $O_\epsilon(n^\epsilon)$. Now apply H\"{o}lder's inequality and get
\begin{align*}
S \ll_{\epsilon, r}& p^{\epsilon / (2r)} \Bigl( \sum_{n \le V H} 1 \Bigr)^{(2r - 1)/(2r)} \Bigl( \sum_{n \le p - 1} \Big| \sum_{u \le U} \Bigl( \frac{u}{p} \Bigr) \Bigl( \frac{u n - 4 c}{p} \Bigr) \Big|^{2r} \Bigr)^{1/(2r)} \\
\ll_{\epsilon, r}& p^{\epsilon / (2r)} (V H)^{1 - 1/(2r)} (U^r p + U^{2r} p^{1/2})^{1/(2r)}
\end{align*}
by Lemma 4 in \cite{Shp} which follows from Weil's bound on multiplicative character sums. Now we take $U = p^{1/(2r)}$ and $V = H / U$ with $1/r < \epsilon$. Then one can verify that $S = o(H^2)$ which implies that there is some $u v \le H / 2$ and $b' \le H / 2$ such that $\Bigl( \frac{u v b'}{p} \Bigr) \Bigl( \frac{u v b' - c}{p} \Bigr) = 1$. This together with the argument in section \ref{basic} gives Theorem \ref{thm1}.

\section{Proof of Theorem \ref{thm2}}
Throughout this section, we assume that $H \gg_\epsilon p^{1/6 + \epsilon}$. Suppose, for all $(c, p) = 1$,
\[
\Bigl( \frac{a' b'}{p} \Bigr) \Bigl( \frac{a' b' - c}{p} \Bigr) = 1
\]
for some $1 \le a', b' \le H/2$. This together with section \ref{basic} implies that
\[
| H_c \cap B(X, Y; H) | \ge 2
\]
for any $(c, p) = 1$ and $H \gg_\epsilon p^{1/6 + \epsilon}$.

\bigskip

Now, suppose this is not the case. Then, for some $(c, p) = 1$,
\[
\Bigl( \frac{a' b'}{p} \Bigr) \Bigl( \frac{a' b' - c}{p} \Bigr) = 0 \text{ or } -1
\]
for all $1 \le a', b' \le H / 2$. Suppose two such pairs give
\[
\Bigl( \frac{a_1' b_1' - c}{p} \Bigr) = 0 = \Bigl( \frac{a_2' b_2' - c}{p} \Bigr).
\]
This implies $| H_c \cap B(X, Y; H) | \ge 2$ automatically by section \ref{basic}. Subsequently, we assume that
\[
\Bigl( \frac{a' b'}{p} \Bigr) \Bigl( \frac{a' b' - c}{p} \Bigr) = -1
\]
for all but at most one pair of $1 \le a', b' \le H / 2$.
This implies that
\[
\Bigl( \frac{a' - c \overline{b'}}{p} \Bigr) = - \Bigl( \frac{a'}{p} \Bigr)
\]
for all but at most one pair of $1 \le a', b' \le H / 2$. Consequently,
\[
\sum_{b' \le H/2} \Big| \sum_{a' \le H/2} \Bigl( \frac{a' - c \overline{b'}}{p} \Bigr) \Big|^{2r} \ge ([H/2] - 1) \Big| \sum_{a' \le H} \Bigl( \frac{a'}{p} \Bigr) \Big|^{2r} =: ([H/2] - 1) |\Sigma|^{2r}.
\]
Now we apply Theorem \ref{Shao} with $N_{b'} = - c \overline{b'}$. First we claim that $c \overline{b_1'} - c \overline{b_2'}$ cannot be congruent to some $l \le H$ modulo $p$ for $1 \le b_1' < b_2' \le H/2$. For otherwise
\[
c \overline{b_1'} - c \overline{b_2'} \equiv l \pmod p
\]
for some $1 \le l \le H$. Let $a_1' \equiv c \overline{b_1'} \pmod p$ and $a_2' \equiv c \overline{b_2'} \pmod p$. Then $(a_1', b_1')$ and $(a_2', b_2')$ would be two points of the modular hyperbola $H_c$ lying in a square of side length $H$ which contradicts our assumption that no such square contains two such points. Therefore, we can apply Theorem \ref{Shao} and get
\[
([H/2] - 1) |\Sigma|^{2r} = \sum_{b' \le H/2} \Big| \sum_{a' \le H/2} \Bigl( \frac{a' - c \overline{b'}}{p} \Bigr) \Big|^{2r} \ll_{\epsilon, r} H^{2r - 2} p^{1/2 + 1/(2r) + \epsilon}.
\]
This implies that
\[
\Sigma \ll_{\epsilon, r} H^{1 - 3/(2r)} p^{(r+1)/(4r^2) + \epsilon / (2r)} = o(H)
\]
if $r$ is sufficiently large. Hence we have that the character sum
\[
\sum_{a' \le H/2} \Bigl( \frac{a'}{p} \Bigr) =o(H).
\]
Feeding this character sum estimate into the standard Vinogradov's trick in obtaining upper bound for the least quadratic nonresidue, we have the first half of Theorem \ref{thm2}.


Tsz Ho Chan \\
Department of Mathematical Sciences \\
University of Memphis \\
Memphis, TN 38152 \\
U.S.A. \\
thchan6174@gmail.com \\

\end{document}